\documentclass[12pt]{amsart}

\usepackage[latin1]{inputenc}
\usepackage{amsmath,amssymb}
\usepackage{latexsym}
\usepackage[top=2in, bottom=1.5in, left=1in, right=1in]{geometry}

\begin{document}

\title{On the Bernstein's constant in convex approximation}

\author{Sorin G. Gal}

\address{Department of Mathematics and Computer Science, University of Oradea,
410087 Oradea, Romania}
\email{galso@uoradea.ro}



\bigskip \bigskip

\begin{abstract}
Denoting by $E_{n}^{(+2)}(f)$ the best uniform approximation of $f$ by convex polynomials of degree $\le n$, there is an open question if there exists the limit $\lim_{n\to \infty}n^{\lambda}E_{n}^{(+2)}(|x|^{\lambda})$ for $\lambda \ge 1$.

\bigskip



\end{abstract}

\maketitle

\section{Introduction}

A famous result of Bernstein \cite{Bern7}, \cite{Bern8} states that
for any $\lambda >0$, $\lambda$ not even integer, there exists finite the limit $\lim_{n\to \infty}n^{\lambda}E_{n}(|x|^{\lambda}) > 0$, where
$$E_{n}(|x|^{\lambda})=\inf\{\max\{|P(x)-|x|^{\lambda}| ; x\in [-1, 1]\} ; P \in {\mathcal{P}}_{n}(\mathbb{R})\},$$
and ${\mathcal{P}}_{n}(\mathbb{R})$ denotes the set of all real polynomials (that is with real coefficients) of degree $\le n$.

It is worth noting that the exact values of the above limits are not known.
Details and generalizations of these results, can be found in the papers \cite{Ganz2}-\cite{Ganz4}, \cite{Lub} and in the book \cite{Ganz1}.

In this note we consider a similar problem in the case of the best approximation by convex polynomials, which was first mentioned in \cite{Gal1}, p. 325.

\section{Convex Approximation}

Taking into account that $|x|^{\lambda}$ is
convex for $\lambda \ge 1$, denoting
$$E_{n}^{(+2)}(f)=\inf \{\|f-P\|_{\infty} ; P\in {\mathcal{P}}_{n}(\mathbb{R}),  P^{\prime \prime}(x)\ge 0, \forall x\in [-1, 1]
\},$$
it is natural to consider the following open question.

{\bf Open Question 1.} {\it There exists finite the limit
$\lim_{n\to \infty}n^{\lambda}E_{n}^{(+2)}(|x|^{\lambda})$ for $\lambda \ge 1$  ?}

Note that since by \cite{KLS1}, for  $f$ convex
on $[-1, 1]$ and $\lambda >0$ we have
$$E_{n}(f)={\mathcal{O}}(n^{-\lambda}) \quad \mbox{ iff } E_{n}^{(+ 2)}(f)={\mathcal{O}}(n^{-\lambda}), \quad n\to \infty,$$
it easily follows that the sequence $(n^{\lambda}E_{n}^{(+
2)}(|x|^{\lambda}))_{n\in \mathbb{N}}$ for $\lambda \ge 1$ is bounded.

More general, we can consider the following.

{\bf Open Question 2.} {\it If $f$ is convex on $\mathbb{R}$ and
the sequence $(\lambda_{n})_{n}$ satisfies the conditions
in \cite{Ganz1},  p. 3, then it is valid that we have
$$\lim_{n\to \infty}E_{n}^{(+2)}\left (f ; L_{\infty}\left [-\frac{n(1-\lambda_{n})}{\sigma},
\frac{n(1-\lambda_{n})}{\sigma}\right ]\right) < \infty ( \mbox{ or
} \le A_{\sigma}^{(+2)}(f, L_{\infty}(\mathbb{R}))<\infty ) \quad
\quad ? $$}
Here $A_{\sigma}^{(+2)}(f, L_{\infty}(\mathbb{R})=\inf \{\|f-g\|_{L_{\infty}(\mathbb{R})} ; g\in B_{\sigma}, g \mbox{ convex on } \mathbb{R}
\}$ and $B_{\sigma}$ denotes the class of all entire functions of exponential type $\sigma$.

This open question is partially supported by the fact that by \cite{LS1}, for a convex function $f$ we have $E_{n}^{(+2)}(f)_{L_{p}}\le E_{n-2}(f^{\prime \prime})_
{L_{p}}$, for all $n\ge 2$ and $p = \infty$ and this
estimate cannot be improved.

\end{document}